\documentclass[a4paper,12pt]{amsart}
\usepackage{amssymb}
\usepackage{ifthen}
\usepackage{cite}
 \usepackage[dvips]{graphicx}
\nonstopmode \numberwithin{equation}{section}
\usepackage{amssymb}
\usepackage{ifthen}
\usepackage{graphicx}
\usepackage{amsmath}
\usepackage[T1]{fontenc} 
\usepackage[utf8]{inputenc}
\usepackage[usenames,dvipsnames]{color}
\usepackage{color}
\usepackage[english]{babel}
\usepackage{fancyhdr}
\usepackage{fancybox}
\usepackage{tikz}
\nonstopmode \numberwithin{equation}{section}
\theoremstyle{plain}

\newtheorem{conj}{Conjecture}

\theoremstyle{definition}
\newtheorem{defi}{Definition}[section]

\newtheorem{thm}{Theorem}[section]
\newtheorem{prob}{Problem}[section]
\newtheorem{cor}{Corollary}[section]
\newtheorem{ques}{Question}[section]
\newtheorem{prop}{Proposition}[section]
\newtheorem{rem}{Remark}[section]
\newtheorem{lem}{Lemma}[section]

\newcounter{minutes}
\divide\time by 60
\newcounter{hours}
\multiply\time by 60
\addtocounter{minutes}{-\time}

\newcounter {own}
\def\theown {\thesection       .\arabic{own}}

\newenvironment{pf}[1][]{%
 \vskip 3mm
 \noindent
 \ifthenelse{\equal{#1}{}}%
  {{\slshape Proof. }}%
  {{\slshape #1.} }%
 }%
{\qed\bigskip}

\newcounter{alphabet}

\def\be{\begin{equation}}
\def\ee{\end{equation}}

\newcommand{\bee}{\begin{enumerate}}
\newcommand{\eee}{\end{enumerate}}

\newcommand{\blem}{\begin{lem}}
\newcommand{\elem}{\end{lem}}
\newcommand{\bthm}{\begin{thm}}
\newcommand{\ethm}{\end{thm}}
\newcommand{\bcor}{\begin{cor}}
\newcommand{\ecor}{\end{cor}}
\newcommand{\beg}{\begin{examp}}
\newcommand{\eeg}{\end{examp}}
\newcommand{\begs}{\begin{examples}}
\newcommand{\eegs}{\end{examples}}

\newcommand{\bdefn}{\begin{defn}}
\newcommand{\edefn}{\end{defn}}

\newcommand{\bprob}{\begin{prob}}
\newcommand{\eprob}{\end{prob}}
\newcommand{\bei}{\begin{itemize}}
\newcommand{\eei}{\end{itemize}}

\newcommand{\bcon}{\begin{conj}}
\newcommand{\econ}{\end{conj}}
\newcommand{\bcons}{\begin{conjs}}
\newcommand{\econs}{\end{conjs}}
\newcommand{\bprop}{\begin{prop}}
\newcommand{\eprop}{\end{prop}}
\newcommand{\br}{\begin{rem}}
\newcommand{\er}{\end{rem}}
\newcommand{\brs}{\begin{rems}}
\newcommand{\ers}{\end{rems}}
\newcommand{\bo}{\begin{obser}}
\newcommand{\eo}{\end{obser}}
\newcommand{\bos}{\begin{obsers}}
\newcommand{\eos}{\end{obsers}}
\newcommand{\bpf}{\begin{pf}}
\newcommand{\epf}{\end{pf}}
\newcommand{\ba}{\begin{array}}
\newcommand{\ea}{\end{array}}
\newcommand{\beq}{\begin{eqnarray}}
\newcommand{\beqq}{\begin{eqnarray*}}
\newcommand{\eeq}{\end{eqnarray}}
\newcommand{\eeqq}{\end{eqnarray*}}

\begin{document}

\title{Operator-valued analogues of multidimensional Bohr-Rogosinski inequalities}

\author{Sabir Ahammed}
\address{Sabir Ahammed, Department of Mathematics, Jadavpur University, Kolkata-700032, West Bengal, India.}
\email{sabira.math.rs@jadavpuruniversity.in}

\author{Molla Basir Ahamed}
\address{Molla Basir Ahamed, Department of Mathematics, Jadavpur University, Kolkata-700032, West Bengal, India.}
\email{mbahamed.math@jadavpuruniversity.in}


\def\thefootnote{}
\footnotetext{ {\tiny File:~\jobname.tex,
printed: \number\year-\number\month-\number\day,
          \thehours.\ifnum\theminutes<10{0}\fi\theminutes }
} \makeatletter\def\thefootnote{\@arabic\c@footnote}\makeatother
\begin{abstract} 
In this article, we first establish operator-valued analogues of multidimensional refined Bohr inequality.
Then we establish operator-valued analogues of multidimensional improved Bohr inequality with a certain power of the norm of the initial coefficient.  Finally, we establish operator-valued analogues of a multidimensional sharp version of Bohr inequality with the initial coefficient being replaced by the norm value of the function. In addition, we establish operator-valued analogues of multidimensional improved Bohr inequality using the quantity $S_r/\pi$. All the results prove to be sharp. 
\end{abstract}
\maketitle
\pagestyle{myheadings}
\markboth{S. Ahammed and M. B. Ahamed}{Operator-valued analogues of multidimensional Bohr-Rogosinski inequalities}
\section{Introduction}
The classical theorem of Bohr \cite{Bohr-1914}, examined a century ago, generates intensive research activity-what is called Bohr phenomenon. Let $\mathcal{B}(\mathcal{H})$ is a Banach algebra of all bounded linear operators on a complex Hilbert spaces $\mathcal{H} $ with the norm $$||A||= \sup_{h\in\mathcal{H}\setminus\{0\}}\dfrac{||Ah||}{||h||}=\sup_{h\in\mathcal{H},||h||=1}||Ah||.$$ Let $ \mathcal{H}^{\infty}(\mathbb{D}, X) $ be the space of bounded analytic functions $f(z)=\sum_{n=0}^{\infty}A_nz^n$ with $||f||_{{\mathbb{H}}^{\infty}(\mathbb{D},X)}:=\sup_{|z|<1}||f(z)||,$ where $A_n\in X$ for all $ n\in\mathbb{N}_0:=\mathbb{N}\cup\{0\} $ from the unit disk $ \mathbb{D}:=\{z\in\mathbb{C} : |z|<1\} $ into a complex Banach space$ X $ and  $ \mathcal{B}^{\prime}(\mathbb{D}, X) $ be the class of functions $f\in\mathcal{H}^{\infty}(\mathbb{D}, X) $ with $||f||_{{\mathbb{H}}^{\infty}(\mathbb{D},X)}\leq1.$  The Bohr radius $ R(X) $ for the class $ \mathcal{B}^{\prime}(\mathbb{D},X) $ is defined by (see \cite{Blasco-2010})
\begin{align*}
	R(X):=\sup\bigg\{r\in (0, 1) : \sum_{n=0}^{\infty}||a_k||r^k\leq 1\; \mbox{for all}\; f(z)=\sum_{n=0}^{\infty}a_kz^k\in\mathcal{B}^{\prime}(\mathbb{D},X), z\in\mathbb{D}\bigg\},
\end{align*}
where $ \sum_{n=0}^{\infty}||a_k||r^k $ is the associated majorant series of function $ f\in\mathcal{H}^{\infty}(\mathbb{D},X) $. The remarkable result of Harald Bohr \cite{Bohr-1914} states that $ R(X)=1/3 $ when $ X=\mathbb{C} $, where norm of $ X $ is the usual modulus of complex numbers. \vspace{1.2mm}

One of our main aims in this paper is to study operator-valued analogues of multidimensional Bohr-Rogosinski  inequalities and their sharpness. In this context, first of all we fix some notations. For the operator $ T\in \mathcal{B}(\mathcal{H})$, we denote $ ||T|| $ as the norm of $ T $ and the corresponding adjoint operator $ T^{*} : \mathcal{H}\rightarrow\mathcal{H} $ of $ T $ is defined by the inner product $ \langle Tx, y\rangle=\langle x, T^{*}y \rangle $ for all $ x, y\in\mathcal{H} $. The operator $ T $ is said to be normal if $ T^*T=T^*T $, self-adjoint if $ T^*=T $, and positive if $ \langle Tx, x \rangle\geq 0 $ for all $ x\in \mathcal{H} $. The identity operator in $ \mathcal{H} $ is denoted by $ I $. 
\vspace{1.2mm} 
\subsection{The classical Bohr inequality and its recent developments}
An interesting classical result related to the family $ \mathcal{B}:=\mathcal{B}^{\prime}(\mathbb{D},\mathbb{D}) $ was discovered by Harald Bohr in 1914 (see \cite{Bohr-1914}) as follows.

\begin{thm}\cite{Bohr-1914}\label{th-1.1}
	If $ f(z)=\sum_{n=0}^{\infty}a_nz^n\in\mathcal{B} $, then $M_f(r):=\sum_{n=0}^{\infty}|a_n|r^n\leq 1 $ for $|z|=r\leq1/3.$
\end{thm}
The number $ K_1 $ defined as the best radius for which this happens, that is, 
$$	K_1:=\sup\left\{r\in (0, 1) : M_f(r)\leq 1,\; \mbox{holds for all}\; f\in\mathcal{B}\;\mbox{with}\; |z|\leq r\right\},$$
is called the Bohr radius $ K_1 $ for $ \mathbb{D} $. Then $ K_1\geq 1/6 $. Subsequently, later M. Riesz, I. Schur and F. Wiener independently improved the inequality $ M_f(r)\leq 1 $ for $|z|\leq1/3$ and found the exact value $ K_1=1/3 $. Henceforth, if there exists a positive real number $ r_0 $ such that the inequality $ M_f(r)\leq 1 $ holds for every element of a class $ \mathcal{F} $ for $ 0\leq r\leq r_0 $ and fails when $ r>r_0 $, then we shall say that $ r_0 $ is an sharp bound for the inequality w.r.t. the class $ \mathcal{F} $. Several other proofs of this interesting inequality were given in different articles  (see \cite{Paulsen-PLMS-2002,Tomic-1962,Sidon-1927}). A detailed account of the development of the topic can be read in the survey article \cite{Ali-Abu-Ponn-2016} and the references therein. In the majorant series $ M_f(r)$ of function $ f\in\mathcal{B} $, the beginning terms play some crucial role in study of the Bohr-type inequalities. For instance, in case of $|a_0|=0$, Tomic \cite{Tomic-1962} has proved  the inequality $ M_f(r)\leq 1 $ holds for $0\leq r \leq {1}/{2}$ and if term $|a_0|$ is replaced by  $|a_0|^2$, then the constant $1/3$ could be replaced by $1/2$. In addition, if  $|a_0|$ is replaced by $|f(z)|$, then the constant $1/3$ could be replaced by $\sqrt{5}-2$ which is best possible (see \cite{Kayu-Kham-Ponnu-2021-JMAA}). \vspace{1.2mm}

In $ 1995 $, Dixon \cite{Dixon & BLMS & 1995} used Bohr's inequality in connection with the long-standing problem of characterizing Banach algebras satisfying von Neumann inequality and since then the Bohr inequality gets much attention from several researchers, and various results are proved for different classes of functions (see e.g.,  \cite{Abu-MN-2013,Aizenberg-Djakov-PAMS-2000,Allu-BSM-2021,Bayart-2014,Boas-1997,Chen-Liu-Ponnusamy-RM-2023,Defant-JFA-2018,Kumar-PAMS-2023,Kumar-Manna-JMAA-2023,Lata-Singh-PAMS-2022,Paulsen-PLMS-2002,Paulsen-BLMS-2006,Ponnusamy-HJM-2021,Ponnusamy-JMAA-2022,Aha-Aha-CMB-2023,Aha-Aha-CMFT-2023,Aha-Aha-RM-2023}).  In recent years, multidimensional analogues of Bohr inequality were considered by several authors (see \cite{Aizeberg-PLMS-2001,Boas-1997,Boas-2000,Liu-Ponnusamy-PAMS-2021,Aha-Aha-CVEE-2023}). Most noticeably, Paulsen and Singh \cite{Paulsen-PAMS-2004} used positivity methods to obtain operator-valued generalizations of Bohr inequality in the single variable case, whereas Defant \emph{et al.} \cite{Defant-JFA-2018} studied the real version of Bohr phenomenon replacing the $ N $-dimensional torus $ \mathbb{T}^N $ with the $ N $-dimensional Boolean cube $ \{\pm 1\}^N:=\{-1, 1\}^N $. However, it can be noted that not every class of functions has the Bohr phenomenon. For example, B$ \acute{a} $n$ \acute{a} $t$ \acute{a} $au \emph{et al.} \cite{Beneteau-2004} showed  that there is no Bohr phenomenon in the Hardy space $ H^p(\mathbb{D}, X) $, where $ p\in [1, \infty) $. Moreover, it is shown in \cite{Liu-Liu-JMAA-2020} that the Bohr inequality of norm type for holomorphic mappings with lacunary series on the unit polydisc in $ \mathbb{C}^n $ does not hold in general. In fact, it is shown in \cite{Liu-Liu-JMAA-2020} that for the function $ f $ given by
\begin{align*}
	f(u)=\left(u_1\frac{u_1-\frac{1}{\sqrt{2}}}{1-\frac{1}{\sqrt{2}}u_1}, u_2\frac{u_2-\frac{2}{\sqrt{5}}}{1-\frac{2}{\sqrt{5}}u_2}\right)^{\prime},\; z=(u_1, u_2)^{\prime}
\end{align*}
that $ f\in \mathcal{H}(\mathbb{D}^2, \overline{\mathbb{D}^2}) $ and 
\begin{align*}
	&\sum_{s=1}^{\infty}\frac{||D^sf(0)\left(u^s\right)||}{s!}>\frac{1}{\sqrt{2}}\left(\frac{2}{\sqrt{5}}+\frac{1}{\sqrt{2}}\right)>1\; \mbox{for}\; u=\left(\frac{1}{\sqrt{2}}, \frac{1}{\sqrt{2}}\right).
\end{align*}
Subsequently, Lin \emph{et al.} \cite{Lin-Liu-Ponnusamy-AMS-2023} have expanded the Bohr inequality by utilizing the Fr$ \acute{u} $chet derivative with a lacunary series in the context of higher-dimensional spaces, specifically, mappings from $ \mathbb{D}^n $ to $ \mathbb{D} $ (or $ \mathbb{D}^n $ to $ \mathbb{D}^n $). Furthermore, in \cite{Lin-Liu-Ponnusamy-AMS-2023}, the authors  have addressed the presence of a constant term in $ f $ and derive two redefined Bohr inequalities in $ \mathbb{D}^n $. This study gives rise to a natural question as follows:
\begin{ques}\label{Q-1.1}
Can we obtain the Bohr phenomenon for holomorphic mappings $ f $ with lacunary series,which take values in a higher dimensional unit ball without some restricted conditions on the mappings $ f $?
\end{ques}
Hamada \cite{Hamada-MN-2023} gave a positive answer the Question \ref{Q-1.1} in the case $ f $ is a holomorphic mapping with lacunary series from the unit ball
$ \mathbb{B} $ of a complex Banach space into the unit ball $ \mathbb{B}_H $ of a complex Hilbert space $ H $ with $ \dim H\geq 2 $. Surprisingly, the Bohr
phenomenon, which we have obtained in the case $ \dim H\geq 2 $, is completely different from that in the case with values in
$ \mathbb{C} $. Also, our result is sharp $ H=l_2 $. Hence, we observe that the Bohr inequality has evolved into a substantial research topic, as it is being explored for various function spaces from different perspectives. The primary focus of this paper is to delve into the operator-valued analogue of multidimensional Bohr-Rogosinski inequalities and their sharpness. Going forward, we provide a concise overview of the research developments concerning the Bohr-Rogosinski inequality and its relationship with the Bohr inequality.
\subsection{Bohr-Rogosinski inequality and Refinements of the classical Bohr inequality}
We now concentrate another well-known result by Rogosinski \cite{Rogosinski-1923}, which states that for any holomorphic self mapping $ f(z)=\sum_{n=0}^{\infty}a_nz^n $ of the open unit disk $ \mathbb{D} $, $ R_f(r):=|\sum_{n=0}^{N-1}a_nz^n|=|f(z)-\sum_{n=N}^{\infty}a_nz^n|\leq 1 $ for $ |z|=r\leq 1/2 $ and for all integer $ N\geq 1 $. The quantity $ 1/2 $ here is best possible and known as the Rogosinski radius and the corresponding inequality $ R_f(r)\leq 1 $ is known as the Rogosinski inequality. Further development in this problem has been made in \cite{Aizenberg-AMP-2012,Aizenberg-SM-2005,Aigenber-CMFT-2009}. In recent years, many questions (including its improved and refined forms) on a related new concept Bohr-Rogosinski phenomenon are currently being studied extensively (see e.g. \cite{Aha-CMFT-2022,Ahamed-AMP-2021,Ahamed-AASFM-2022,Das-JMAA-2022,Huang-Liu-Ponnu-AMP-2020,Liu-Liu-Ponnusamy-2021}). However, to the best of our knowledge, unlike the refined Bohr inequalities, there is no multidimensional analogue version of the Rogosinski inequality.  Motivated by this, in this paper, our prime concern is to study the multidimensional analogue of the Bohr-Rogosinski inequality and its sharpness. \vspace{1.2mm}
We define 
\begin{align*}
	\mathcal{T}^N_{j, f}(r):=|f(z)|^j+\sum_{n=N}^{\infty}|a_n|r^n\; \mbox{for}\; j=1, 2.
\end{align*}
Inspired by the notion of the Rogosinski inequality and the Rogosinski’s radius investigated in \cite{Rogosinski-1923}, in \cite{Kayu-Kham-Ponnu-2021-JMAA} author's have obtained the following Bohr-Rogosinski inequality and Bohr-Rogosinski radius for the class $ \mathcal{B} $.
\begin{thm}\cite{Kayu-Kham-Ponnu-2021-JMAA}\label{th-11.99}
	If $ f(z)=\sum_{n=0}^{\infty}a_nz^n\in\mathcal{B} $, then $ 	\mathcal{T}^N_{1, f}(r)\leq 1\;\;\mbox{for}\;\; r\leq R_N, $	where $ R_N $ is the positive root of the equation $ 2(1+r)r^N-(1-r)^2=0 $. The radius $ R_N $ is the best possible. Moreover, $ \mathcal{T}^N_{2, f}(r)\leq 1\;\;\mbox{for}\;\; r\leq R^{\prime}_N, $	where $ R^{\prime}_N $ is the positive root of the equation $ (1+r)r^N-(1-r)^2=0 $. The radius $ R^{\prime}_N $ is the best possible.
\end{thm}
In \cite{Ism-Kay-Ponn-JMAA-2020}, Ismagilov \email{et al.} observed that $ \mathcal{T}^1_{1, f}(r)\leq 1 $ for $ r\leq \sqrt{5}-2 $ and $ \mathcal{T}^2_{2, f}(r)\leq 1 $ for $ r\leq 1/3 $, where the constants $ \sqrt{5}-2 $ and $ 1/3 $ are best possible.\vspace{1.2mm}

In a number of articles (e.g., \cite{Das-JMAA-2022,Huang-Liu-Ponnu-AMP-2020,Liu-Liu-Ponnusamy-2021}), Theorem \ref{th-11.99} is refined further and sometimes improved it for different classes of functions  establishing sharp results. Research on finding sharp results continues in view of refining the Bohr inequality for the class $ \mathcal{B} $ of self-analytic maps on the unit disk $ \mathbb{D} $ and now-a-days researchers are mainly engaged in finding possible sharp refined Bohr inequality. For $ p=1, 2 $, some of the refined versions of L.H.S of the classical Bohr inequality can be defined as the following
\[
\begin{cases}
	\mathcal{A}_{p, f}(r):=|f(z)|^{p}+\displaystyle\sum_{n=1}^{\infty}|a_n|r^n+\left(\dfrac{1}{1+|a_0|}+\dfrac{r}{1-r}\right)\sum_{n=1}^{\infty}|a_n|^2r^{2n}\vspace{1.2mm}\\
	\mathcal{B}_{p, f}(r):=\displaystyle\sum_{n=0}^{\infty}|a_n|r^n+\left(\dfrac{1}{1+|a_0|}+\dfrac{r}{1-r}\right)\sum_{n=1}^{\infty}|a_n|^2r^{2n}+|f(z)-a_0|^{p}.
\end{cases}
\]
We recall here such results as the refined formulation of the Bohr-Rogosinski inequality which are established by Liu \emph{et al.} \cite{Liu-Liu-Ponnusamy-2021}. In fact, the authors studied several sharp refinements of the classical Bohr inequality with various settings. Hereafter, we revisit the results that we intend to examine within the operator-valued setting.
\begin{thm}\label{th-1.5} 
	\cite{Liu-Liu-Ponnusamy-2021} Suppose that $f=\sum_{n=0}^{\infty}a_nz^n\in\mathcal{B}.$
	\begin{enumerate}
		\item[(i)] Then $ \mathcal{A}_{1, f}(r)\leq 1 $ for $|z|=r \leq r_0=2/\left(3+|a_0|+\sqrt{5}(1+|a_0|)\right).$ The radius $r_0$ is best possible and $r_0>\sqrt{5}-2.$ Moreover, $ \mathcal{A}_{2, f}(r)\leq 1 $ for $|z|=r\leq r^{\prime}_0,$ where $r^{\prime}_0$ is the unique positive root of the equation 
		$(1-|a_0|^3)r^3-(1+2|a_0|)r^2-2r+1=0.$ The radius $r^{\prime}_0$ is best possible. Further, we have $1/3<r^{\prime}_0<1/\left(2+|a_0|\right). $\vspace{1mm}
		\item[(ii)] In addition, $ \mathcal{B}_{1, f}(r)\leq 1 $ for $ r\leq 1/5 $, and $ \mathcal{B}_{2, f}(r)\leq 1 $ for $ r\leq 1/3 $ if and only if $ 0\leq |a_0|\leq 4\sqrt{2}-5 $, where the constants $ 1/5 $ and $ 1/3 $ cannot be improved.
	\end{enumerate}
\end{thm} 
In this article, we continue the study to give operator-valued analogues of multidimensional Bohr-Rogosinski inequalities and their sharp refinements and we discuss this part in the next section.
\section{{Operator-valued analogues of multidimensional Bohr-Rogosinski inequality}}
Unlike the study of Bohr inequalities for certain classes of analytic functions of one complex variable, the same is studied by several researchers in the case of functions of $ n $-complex variables. Before we proceed further, we need to introduce some concepts. Here and hereafter, we use the standard multi-index notation  $\alpha$ which is an $n$-tuple $(\alpha_1,\alpha_2,\dots \alpha_n)$ of non-negative integers, $|\alpha|$ is the sum $\alpha_1+\alpha_2+\dots+\alpha_n$ of its components, $\alpha!$ is the product $\alpha_1!\alpha_2!\dots \alpha_n!,$  $z$ denotes an $n$-tuple $(z_1,z_2,\dots z_n)$ of complex numbers, and $z^{\alpha}$ denotes the product $z^{\alpha_1}_1z^{\alpha_2}_2\dots z^{\alpha_n}_n. $
Let $\mathbb{D}^n:=\{z\in\mathbb{C}^n:z=(z_1,z_2, \dots ,z_n), |z_j|<1,j=1,2,\dots,n\}$ be the open unit polydisk. Also, let ${\mathcal{K}_n}$ be the largest non negative number such that if the $n$- variables power series $\sum_{\alpha}a_{\alpha}z^{\alpha}$ 
converges in $\mathbb{D}^n$ and its sum $f$ has modulus less than $1$ so that $a_{\alpha}={\partial}^{\alpha}f(0)/\alpha!,$ then $ \sum_{\alpha}|a_\alpha||z^\alpha|<1\;\;\mbox{for all}\;\; z\in {\mathcal{K}_n}.\mathbb{D}^n. $
\begin{defi}
	A domain $D\subset\mathbb{C}^n$ is said to be a Reinhardt domain centered at $0\in D$ if for any $z=(z_1,z_2,\dots,z_n)\in D,$ we have that $(z_1e^{i\theta_1},z_2e^{i\theta_2},\dots,z_ne^{i\theta_n})\in D$ for each $\theta_k\in[0,2\pi], k=1,2,\dots,n$. We say that $D\subset\mathbb{C}^n$ is a complete Reinhardt domain if for each $z=(z_1,z_2,\dots,z_n)\in D,$ and for each $|u_k|\leq1,$ $k=1,2,\dots,n,$ we have that $u.z=(u_1z_1,u_2z_2,\dots,u_nz_n)\in D.$
\end{defi}
\begin{defi}
	A domain $Q\subset \mathbb{C}^n$ is called a circular domain centered at $0\in Q$ if for any $z=(z_1,z_2,\dots,z_n)\in Q,$ and for each $\theta\in [0,2\pi],$ we have that $(z_1e^{i\theta},z_2e^{i\theta},\dots,z_ne^{i\theta})\in Q.$ We say that  $Q\subset \mathbb{C}^n$ is a complete circular domain centered at $0\in Q$ if  for each $z=(z_1,z_2,\dots,z_n)\in Q,$ and for each $|u|\leq1,$  we have that $u.z=(uz_1,uz_2,\dots,uz_n)\in Q.$
\end{defi}
If $Q\subset{\mathbb{C}^n}$ be a complete circular domain  centered at $0\in Q ,$  every holomorphic function $f$ in $Q$ can be expanded into homogeneous polynomials given by 
\begin{align}\label{e-2.1}
	f(z)=\sum_{k=0}^{\infty}P_k(z)\;\;\mbox{for }\;\; z\in Q,
\end{align}
where $P_k(z)$ is a homogeneous polynomial of degree $k,$ and $P_0(z)=f(0).$ \vspace{1.2mm}

\par In 1997, Boas and Khavinson \cite{Boas-1997} introduced the $ N $-dimensional Bohr radius $ K_N, $ $ (N\in\mathbb{N}$ ) for the polydisk $ \mathbb{D}^N=\mathbb{D}\times\cdots\times\mathbb{D} $. In fact, Boas and Khavinson \cite{Boas-1997} proved that the $ N $-Bohr radius $ K_N $ as the largest radius $ r>0 $ such that for every complex polynomials $ \sum_{\alpha\in\mathbb{N}_0^N}c_{\alpha}z^{\alpha} $ in $ N $ variables 
\begin{align*}
	\sup_{z\in r\mathbb{D}^N} \sum_{\alpha\in\mathbb{N}_0^N}|c_{\alpha}z^{\alpha}|\leq \sup_{z\in \mathbb{D}^N} \bigg|\sum_{\alpha\in\mathbb{N}_0^N}c_{\alpha}z^{\alpha}\bigg|.
\end{align*}
As expected, the constant $ K_N $ is defined as the largest radius $ r $ satisfying $ \sum_{\alpha}|c_{\alpha}z^{\alpha}|<1 $ for all $ z $ with $ ||z||_{\infty}:=\max\{|z_1|, |z_2|, \ldots, |z_N|\}<r $ and all $ f(z)=\sum_{\alpha}c_{\alpha}z^{\alpha}\in\mathcal{H}(\mathbb{D}^N) $. This paper stimulated interests on Bohr-type questions in multidimensional analogues of the Bohr-type inequalities. The reader is refereed to  (\cite{Aizn-OTAA-2005,Aizn-PAMS-2000,Aizn-ST-2007,Aizenberg-Djakov-PAMS-2000,ADLF-IJM-2006,Liu-Ponnusamy-PAMS-2021}) and references therein for results with different settings of the classical Bohr inequality.\vspace{1.2mm}

The recent development of the Bohr inequality is  established as an operator-valued analogue of the classical Bohr inequality and this study opens a new dimension in further investigations of the Bohr-type inequalities. For instance, Popescu \cite{Popescu-2019},  Allu and Halder \cite{{Allu-CMB-2022}} have studied the operator-valued Bohr phenomenon in view of certain general settings. More precisely, Popescu \cite{Popescu-2019} has shown that the classical theorem of Bohr can be proved analogously for operator-valued bounded analytic functions in $\mathbb{D}.$
\begin{thm}\label{thh-3.1}\cite{Popescu-2019}
	Let $f(z)=\sum_{n=0}^{\infty}A_nz^n\in{{\mathcal{H}}^{\infty}(\mathbb{D},\mathcal{B}(\mathcal{H}))}$ such that $A_0=a_0I,$ $a_0\in\mathbb{D}.$ Then
	\begin{align}\label{ee-3.1}
		\sum_{n=0}^{\infty}||A_n||r^n\leq{||f||}_{{\mathbb{H}}^{\infty}(\mathbb{D},\mathcal{B}(\mathcal{H}))}\;\;\mbox{for}\;\; |z|=r\leq\dfrac{1}{3},
	\end{align}
	and $1/3$ is the best possible constant. Moreover, the inequality is strict unless $f$ is a constant.
\end{thm} 
Using the standard multi-index notation, we write an operator-valued $n$-variable power series $ f(z)=\sum_{\alpha}A_\alpha z^\alpha, \;\; A_\alpha\in\mathcal{B}(\mathcal{H}). $ Let $\mathcal{K}_n(\mathcal{H})$ be the largest non-negative number such that this power series  converges in $\mathbb{D}^n$ and ${||f||}_{{\mathcal{H}}^{\infty}(\mathbb{D},\mathcal{B}(\mathcal{H}))}\leq1,$ then $ {\sum_{\alpha}||A_\alpha||\; |z^\alpha|}\leq 1\;\;\mbox{for all}\;\; z\in\mathcal{K}_n(\mathcal{H}).\mathbb{D}^n. $ 
\begin{defi}
	The class $ {\mathcal{B}}(\mathbb{D},\mathcal{B}(\mathcal{H})) $ is said to satisfy the Bohr phenomenon if all the functions in $ {\mathcal{B}}(\mathbb{D},\mathcal{B}(\mathcal{H})) $ satisfies the inequality \eqref{ee-3.1} for $ r\leq 1/3 $. It is worth mentioning that the constant $ 1/3 $ independent of the coefficients of the functions.
\end{defi}
With the help of Counterpart of Schwarz-Pick inequality in \cite[Lemma 1.6]{AR-2006, Allu-CMB-2022}, it is shown that if $f(z)=\sum_{n=0}^{\infty}A_nz^n\in{{\mathcal{B}}(\mathbb{D},\mathcal{B}(\mathcal{H}))}$ such that $A_0=a_0I,$ where  $a_0\in\mathbb{D},$ then 
\begin{align}\label{ee-3.5}
	||A_n||\leq||I-|A_0|^2||=1-|a_0|^2\;\;\mbox{for all}\;\;  n\in\mathbb{N}.
\end{align} 
In view of the above estimate, Allu and Halder \cite{Allu-CMB-2022} continued further study on the result of Popescu \cite{Popescu-2019} and proved the following result.
\begin{thm}\cite{Allu-CMB-2022} \label{th-2.22}
	If the series \eqref{e-2.1} converges in the domain $ Q $ such that the estimate $ ||f(z)||<1 $ holds in $ Q $ and $ f(0)=a_0I $, where $ a_0\in\mathbb{C} $, $ |a_0|<1 $, then $\sum_{k=0}^{\infty} ||P_k(z)||<1 $ in the homothetic domain $ (1/3)Q $. Moreover, if $ Q $ is convex, then $ 1/3 $ is best possible.
\end{thm}
In \cite[Lemma 3.3]{Allu-CMB-2022}, Allu and Halder  showed that if $f(z)=\sum_{n=0}^{\infty}A_nz^n\in{{\mathcal{H}}(\mathbb{D},\mathcal{B}(\mathcal{H}))}$ such that $||f(0)||<1,$ then
\begin{align}\label{ee-3.6}
	||f(z)||\leq\dfrac{||f(0)||+|z|}{1+||f(0)||\;|z|}\;\;\mbox{for}\;\;z\in\mathbb{D},
\end{align}
and in view of this estimate, the authors obtained the following sharp results.
\begin{thm}\cite{Allu-CMB-2022}\label{thh-3.2}
	If $Q$ and $ f $ are as in the Theorem  \ref{th-2.22}, then for $p\in(0,1],$ we have $ {||f(z)||}^p+\sum_{k=N}^{\infty}||P_k(z)||\leq 1 $
	in the homothetic domain $(R_{N,p})Q,$ where  $(R_{N,p})$ is the root in $(0,1)$ of the equation $ 2(1+r)r^N-p(1-r)^2=0. $
	Moreover, if $Q$ is convex, then the constant  $(R_{N,p})$ is the best possible. 
\end{thm}
\begin{thm}\cite{Allu-CMB-2022} \label{thh-3.3}
	If $Q$ and $ f $ are as in the Theorem \ref{thh-3.2},  then  
	\begin{align}\label{ee-2.7}
		||f(0)||+\sum_{k=1}^{\infty}{||P_k(z)||}+\left(\dfrac{1}{1+||f(0)||}+\dfrac{r}{1-r}\right)\sum_{k=1}^{\infty}{||P_k(z)||}^2\leq 1
	\end{align} 
	holds in the homothetic domain $(1/3)Q.$ Moreover, if $Q$ is convex, then $1/3$ is the best possible.
\end{thm}
Let $ \mathbb{D}_r:=\{z\in\mathbb{D} : |z|<r\} $ and $S_r=S_f(r)$ be the planar integral $ S_r:=\int_{\mathbb{D}_r}|f^{\prime}(z)|dA(z) $ for the holomorphic function $ f(z)=\sum_{n=0}^{\infty}a_nz^n $ defined on $ \mathbb{D} $. Then it is well-known that $S_r/\pi=\sum_{k=1}^{\infty}k|a_k|^2r^{2k}.$   Several results are established as the improved versions of the Bohr inequality with the help of the quantity $S_r/\pi$ (see e.g.  \cite{Ahamed-AASFM-2022,Huang-Liu-Ponnu-CVEE-2021,Ism-Kay-Ponn-JMAA-2020,Kayumov-CRACAD-2018}). Like the quantity $S_r/\pi,$ we define 
\begin{align}\label{e-2.8}
	\dfrac{S_r}{\pi}:=\sum_{n=1}^{\infty}n{||A_n||}^2r^{2n}, \;\;\mbox{for}\;\; |z|=r<1,
\end{align}  
for the function $f(z)=\sum_{n=0}^{\infty}A_nz^n\in{{\mathcal{B}}(\mathbb{D},\mathcal{B}(\mathcal{H}))}.$ Let $G(t)$ be a polynomial defined  by 
\begin{align}\label{e-2.9}
	G(t)=c_1t+c_2t^2+\dots+c_lt^l, \;\mbox{for}\; c_i\in\mathbb{R}^{+}\cup\{0\}, i=1,2,\dots,l.
\end{align} 
We define $ \mathcal{E}_f(r):=||f(0)||+\sum_{k=1}^{\infty}{||P_k(z)||}+G\left(\sum_{k=1}^{\infty}k{||P_k(z)||}^2\right). $
\begin{thm}\cite{Allu-CMB-2022} \label{thm-2.5}
	If $Q$ and $ f $ are as in the Theorem \ref{th-2.22}, then  $ \mathcal{E}_f(r)\leq 1 $ 
	holds in the homothetic domain $(1/3)Q,$ where $G(t)$ is given by \eqref{e-2.9} and  the coefficients of $G(t)$ satisfy $8c_1\left({3}/{8}\right)^2+24c_2\left({3}/{8}\right)^4+\dots+8(2l-1)c_l\left({3}/{8}\right)^{2l}\leq1$. Moreover, if $Q$ is convex, then $1/3$ cannot be replaced by a larger quantity. 
\end{thm}
In this article, we are interested to consider improved versions of the Bohr-Rogosinski inequality in the multidimensional case. Thus, to continue the study, and also to give an operator-valued analogue of the results holds for the class $ \mathcal{B} $, it is natural to raise the following. 
\begin{prob}\label{p-3.1}
	Can we establish the analogues of Theorem \ref{th-1.5} for operator-valued functions in $ \mathcal{B}(\mathbb{D}, \mathcal{B}(\mathcal{H})) $? If so, then what is the multidimensional analogue of Theorem \ref{th-1.5} for functions in the class $ \mathcal{B}(\mathbb{D}, \mathcal{B}(\mathcal{H})) $?
\end{prob}
\begin{prob}\label{p-3.2}
	Can we establish an improved version of the inequality \eqref{ee-2.7} in Theorem \ref{thh-3.3}?
\end{prob}
In addition, to obtain an improved version of Theorem \ref{thm-2.5}, we are interested to solve the following.
\begin{prob}\label{p-3.3}
	Can we establish an improved version of Theorem \ref{thm-2.5} replacing the quantity $ ||f(0)|| $ by $ {||f(z)||} $?
\end{prob}
Moreover, to obtain an improved version of Theorem \ref{thh-3.3}, we are interested in solving the following.
\begin{prob}\label{p-2.4}
	Can we establish an improved version of Theorem \ref{thh-3.3} using the quantity $S_r/\pi?$	
\end{prob}
Inspired by the results in \cite{Popescu-2019} and \cite{Allu-CMB-2022}, in this paper, we obtain sharp results in solving the above problems completely. The organization of the paper is the following: In Section 3, we state some key lemmas and the main results of the paper and discussed some relevant remarks. In Section 4, we discuss the proof of the key lemmas and the main results. 

\section{Key lemmas and main results}
In this section, we first state some lemmas which will play a key role to prove the corresponding main results of this paper. Henceforth, for $ j=1, 2 $, we introduce some notations as follows:
\[
\begin{cases}
	\mathcal{M}_{j, f}(r):=\displaystyle ||f(z)||^{j}+\sum_{n=1}^{\infty}||A_n||r^n+\left(\dfrac{1}{1+||A_0||}+\dfrac{r}{1-r}\right)\sum_{n=1}^{\infty}{||A_n||}^2r^{2n}\vspace{0.5mm}\\ \displaystyle \mathcal{C}_{j, f}(r):=\sum_{n=0}^{\infty}||A_n||r^n+\left(\dfrac{1}{1+||A_0||}+\dfrac{r}{1-r}\right)\sum_{n=1}^{\infty}{||A_n||}^2r^{2n}+||f(z)-A_0||^j\vspace{0.5mm}\\
	\displaystyle \mathcal{D}_{ f}(r):={||f(z)||}+\sum_{n=1}^{\infty}||A_n||r^n+G\left(\dfrac{S_r}{\pi}\right)\\
	\mathcal{N}_f^i(r):={||A_0||}^j+\sum_{n=1}^{\infty}||A_n||r^n+\left(\dfrac{1}{1+||A_0||}+\dfrac{r}{1-r}\right)\sum_{n=1}^{\infty}{||A_n||}^2r^{2n}+\lambda\left(\dfrac{S_r}{\pi}\right),
\end{cases}
\]
where $S_r/\pi$ is given by \eqref{e-2.8}.\vspace{1.2mm}

We prove the Lemmas \ref{lem-33.11} and \ref{lem-4.3} which are operator-valued analogue of part (i) and (ii), respectively,  of Theorem \ref{th-1.5}. This solves the first part of Problem \ref{p-3.1} completely.
\begin{lem}\label{lem-33.11}
	Suppose that $f\in\mathcal{B}(\mathbb{D},\mathcal{B}(\mathcal{H}))$ with the expansion $f(z)=\sum_{n=0}^{\infty}A_nz^n$ in $\mathbb{D}$ such that $A_n\in\mathcal{B}(\mathcal{H})$ for all  $ n\in\mathbb{N}_0 $ and $f(0)=b_0I,$ $b_0\in\mathbb{C},$ $|b_0|<1.$ If $||f(z)||\leq1$ in $\mathbb{D},$ then $ \mathcal{M}_{1, f}(r)\leq 1 $ for $|z|=r\leq1/\sqrt{5}.$
	The constant $1/\sqrt{5}$ is best possible.\vspace{1.2mm}
	Moreover,	$ \mathcal{M}_{2, f}(r)\leq 1 $ for $|z|=r\leq1/3 .$ The constant $1/3$ is best possible.
\end{lem}
\begin{lem}\label{lem-4.3}
	Suppose that $f$ be as in Lemma \ref{lem-33.11}. If $||f(z)||\leq1$ in $\mathbb{D},$ then  $ \mathcal{C}_{1, f}(r)\leq 1$ for $|z|=r\leq1/5.$	The constant $1/5$ is the best possible. Moreover, $ C_{2,f}(r)\leq 1 $ for $|z|=r\leq1/3.$ The constant $1/3$ is best possible.
\end{lem}
We prove the following two lemmas which we will use to prove results in order to solve the Problems \ref{p-3.2} and \ref{p-3.3}, respectively.
\begin{lem}\label{lem-2.3}
	Suppose that $f$ be as in Lemma \ref{lem-33.11}. If $||f(z)||\leq1$  in $\mathbb{D}$ and $p\in(0,1],$ then 
	\begin{align}\label{ee-4.4}
		\mathcal{B}_p(r):={||A_0||}^p+\sum_{n=1}^{\infty}||A_n||r^n+\left(\dfrac{1}{1+||A_0||}+\dfrac{r}{1-r}\right)\sum_{n=1}^{\infty}{||A_n||}^2r^{2n}\leq 1
	\end{align} for $|z|=r\leq1/3.$
	The constant $1/3$ is the best possible.
\end{lem}
\begin{lem}\label{lem-3.4}
	Suppose that $f$ be as in Lemma \ref{lem-33.11}. If $||f(z)||\leq 1$ in $\mathbb{D},$ then the following sharp inequality holds: $ \mathcal{D}_{f}(r)\leq1\;\;\mbox{for}\;\; r\leq r_0=\sqrt{5}-2, $ where $G(t)$ is a polynomial defined by \eqref{e-2.9} and satisfy 
	\begin{align}\label{eee-3.8}
		\sum_{m=1}^{l}c_m2^{(1-4m)}\leq\dfrac{13-5\sqrt{5}}{4}.
	\end{align}
\end{lem}
The quantity $ S_r/\pi $ has a certain significant role in the study of improved versions of  Bohr-type inequalities (see e.g. \cite{Ahamed-AASFM-2022,Ism-Kay-Ponn-JMAA-2020,Liu-Liu-Ponnusamy-2021,Huang-Liu-Ponnu-CVEE-2021,Kayumov-CRACAD-2018}). To establish operator-valued analogues of the multidimensional improved Bohr inequality using the quantity $S_r/\pi,$ we establish the next two Lemmas.
\begin{lem}\label{BS-lem-3.5}
	Suppose that $f$ be as in Lemma \ref{lem-33.11}. If $||f(z)||\leq 1$ in $\mathbb{D},$ then the following inequality holds: 	
	\begin{align}\label{BS-eq-3.4}
		\dfrac{S_r}{\pi}\leq\frac{r^2(1-||A_0||^2)^2}{(1-||A_0||^2r^2)^2}\;\;\;\;\; \mbox{for}\; r\leq \frac{1}{\sqrt{2}}.
	\end{align}
	where $S_r/\pi$ is the quantity given by \eqref{e-2.8}.
\end{lem}
We establish the following lemma, which will serve as a basis for proving results essential to the resolution of Problem \ref{p-2.4}.
\begin{lem}\label{BS-lem-3.6}
	Suppose that $f$ be as in Lemma \ref{lem-33.11}. If $||f(z)||\leq1$  in $\mathbb{D}$, then $\mathcal{N}_f^1(r)\leq 1$
	for $|z|=r\leq1/3,$ $\lambda=8/9,$ and $0\leq||A_0||\leq0.402964.$
	The constants $1/3$ and $8/9$ are best possible. Moreover, $\mathcal{N}_f^2(r)\leq 1$
	for $ r\leq 1/(3-||A_0||) ,$ $\lambda=9/8$ and $0\leq||A_0||\leq 0.489758.$  The constants $ 1/(3-||A_0||) $ and $ 9/8 $ cannot be improved.	
\end{lem}
Before we state the main results of the paper, we introduce here some notations. For $ j=1, 2 $, we define
\[
\begin{cases}
	\mathcal{I}_{j, f}(r):=\displaystyle||f(z)||^j+\sum_{k=1}^{\infty}{||P_k(z)||}+\left(\dfrac{1}{1+||f(0)||}+\dfrac{r}{1-r}\right)\sum_{k=1}^{\infty}{||P_k(z)||}^2\vspace{0.5mm}\\
	\mathcal{J}_{j, f}(r):=\displaystyle\sum_{k=0}^{\infty}{||P_k(z)||}+\left(\dfrac{1}{1+||f(0)||}+\dfrac{r}{1-r}\right)\sum_{k=1}^{\infty}{||P_k(z)||}^2+||f(z)-f(0)||^j\\
	\mathcal{L}_{ f}(r):=\displaystyle {||f(z)||}+\sum_{k=1}^{\infty}{||P_k(z)||}+G\left(\sum_{k=1}^{\infty}k{||P_k(z)||}^2\right)\\
	\mathcal{R}^j_f(r):={||f(0)||}^j+\displaystyle\sum_{k=1}^{\infty}{||P_k(z)||}+\left(\dfrac{1}{1+||f(0)||}+\dfrac{r}{1-r}\right)\sum_{k=1}^{\infty}{||P_k(z)||}^2\\\quad\quad\quad\quad+\lambda \left(\sum_{k=1}^{\infty}k{||P_k(z)||}^2\right).
\end{cases}
\]
We obtain the following result which is multidimensional analogue of part (i) of Theorem \ref{th-1.5}. 
\begin{thm}\label{thm-3.1} 
	Suppose that $Q$ is a complete circular domain centered at  $0\in Q \subset{\mathbb{C}^n}.$  If  the series $\eqref{e-2.1}$ converges in the domain $Q$  such that $||f(z)||<1$ for all $z\in Q,$ and $||f(0)||<1.$ Then  $ \mathcal{I}_{1, f}(r)\leq 1 $ holds in the homothetic domain $(1/\sqrt{5})Q.$ Moreover, if $Q$ is convex, then $1/\sqrt{5}$ is the best possible. Moreover, $ \mathcal{I}_{2, f}(r)\leq 1 $ holds in the homothetic domain $(1/3)Q.$ Moreover, if $Q$ is convex, then $1/3$ is the best possible.
\end{thm}
We establish the following result  which is multidimensional analogue of part (ii) of Theorem \ref{th-1.5}.
\begin{thm}\label{thm-3.3}
	Suppose that $Q$ and $f$ be as in Theorem \ref{thm-3.1}. Then $ \mathcal{J}_{1, f}(r)\leq 1 $ holds in the homothetic domain $(1/5)Q.$ Moreover, if $Q$ is convex, then $1/5$ is the best possible. Moreover, $ \mathcal{J}_{2, f}(r)\leq 1 $ holds in the homothetic domain $(1/3)Q.$ Moreover, if $Q$ is convex, then $1/3$ is the best possible.
\end{thm}
\begin{rem}
	It is clear from Theorems \ref{thm-3.1} and \ref{thm-3.3} that the second part of Problem \ref{p-3.1} is solved completely.
\end{rem}
We prove the following result which generalizes Theorem \ref{thh-3.3} considering certain powers of the initial coefficients in the majorant series of $ f $ given by \eqref{e-2.1}.
\begin{thm}\label{thm-3.2}
	Suppose that $Q$ and $f$ be as in Theorem \ref{thm-3.1} and $p\in(0,1]$.  Then 
	\begin{align}\label{ee-33.5}
		{||f(0)||}^p+\sum_{k=1}^{\infty}{||P_k(z)||}+\left(\dfrac{1}{1+||f(0)||}+\dfrac{r}{1-r}\right)\sum_{k=1}^{\infty}{||P_k(z)||}^2\leq 1 
	\end{align} 
	holds in the homothetic domain $(1/3)Q.$ Moreover, if $Q$ is convex, then $1/3$ is the best possible. 
\end{thm}
\begin{rem}
	It is easy to see that Theorem \ref{thh-3.3} still sharp if the first term $ ||f(0)|| $ in L.H.S of \eqref{ee-2.7} is replaced by $ ||f(0)||^p $. Clearly, Theorem \ref{thm-3.2} improves the Theorem \ref{thh-3.3} as the Problem \ref{p-3.2} is solved successfully  for the values $ p\in[0,1) $.
\end{rem} 
We obtain the following result solving the Problem \ref{p-3.3}.
\begin{thm}\label{thh-3.4}
	Suppose that $Q$ and $f$ be as in Theorem \ref{thm-3.1}.   If $G(t) $ given by \eqref{e-2.9}, then  $ \mathcal{L}_{ f}(r)\leq 1 $ holds in the homothetic domain $(\sqrt{5}-2)Q,$ where the coefficients of $G(t)$ satisfy \eqref{eee-3.8}. Moreover, if $Q$ is convex, then the constant $\sqrt{5}-2$ is best possible.
\end{thm}
The following is an immediate corollary of Theorem \ref{thh-3.4} which is a multidimensional analogue of \cite[Theorem 3, pp. 4]{Ism-Kay-Ponn-JMAA-2020}.
\begin{cor} Let f be as in Theorem \ref{thh-3.4}. If $G$ is a monomial of one degree,  then  
	\begin{align*}
		{||f(z)||}+\sum_{k=1}^{\infty}{||P_k(z)||}+2(\sqrt{5}-1)\left(\sum_{k=1}^{\infty}k{||P_k(z)||}^2\right)\leq 1 
	\end{align*}
	holds in the homothetic domain $(\sqrt{5}-2)Q.$ Moreover, if $Q$ is convex, then the constant $\sqrt{5}-2$ is the best possible. 
\end{cor}
We obtain the following result which is
improved version of Theorem \ref{thh-3.3} and the Problem \ref{p-2.4} is solved successfully.
\begin{thm}\label{BS-thm-3.5}
	Suppose that $Q$ and $f$ be as in Theorem \ref{thm-3.1}.  Then $\mathcal{R}^1_f(r)\leq 1$
	holds in the homothetic domain $(1/3)Q,$ $\lambda=8/9,$ and $ 0\leq a=||f(0)||\leq 0.402964$. Moreover, if $Q$ is convex, then $1/3$ and  $8/9$ are best possible. In addition, 
	$\mathcal{R}^2_f(r)\leq 1$
	holds in the homothetic domain $(1/(3-||f(0)||))Q,$  $\lambda=9/8,$ and $0\leq ||f(0)||\leq0.489758. $ Moreover, if $Q$ is convex, then $1/(3-||f(0)||)$ and $9/8$ are the best possible.
\end{thm}
\begin{remark}
	In fact,  Theorem \ref{BS-thm-3.5} is multidimensional analogue of \cite[Theorem 4, pp.12]{Liu-Liu-Ponnusamy-2021}.
\end{remark}
\section{Proof of the lemmas and main results}
In this section, we first give a detailed proof of all the lemmas, and later, we discuss proof of the main results of the paper. 
\begin{proof}[\bf Proof of Lemma \ref{lem-33.11}.]
	Let $a=||A_0||.$ 	In view of \eqref{ee-3.5} and \eqref{ee-3.6}, a simple computation shows that
	\begin{align*}
		\mathcal{M}_{1, f}(r)&\leq \left(\dfrac{a+r}{1+ra}\right)+\dfrac{(1-a^2)r}{1-r}+\dfrac{(1+ra)(1-a^2)^2r^2}{(1+a)(1-r)(1-r^2)}\\&=1+\dfrac{(1-a)P_1(a,r)}{(1+ar)(1-r)(1-r^2)},
	\end{align*}
	where $ 	P_1(a,r):=(r-1)^2(1-r^2)+(1+ar)(1+a)r(1-r^2)+(1+ar)^2(1-a^2)r^2. $ We claim that $P_1(a,r)\leq 0$ for $r\leq 1/\sqrt{5}$ and $a\in[0,1).$
	Thus it follows that 
	\[
	\begin{cases}
		\dfrac{\partial P_1(a,r)}{\partial a}=r(1-r^2)(1+r+2ra)+r^2\left(2r+2a(r^2-1)-6ra^2-4a^3r^2\right),\vspace{2.3mm}\\
		\dfrac{{\partial}^2 P_1(a,r)}{\partial a^2}=2r^2(1-r^2)+r^2\left(2(r^2-1)-12ra-12a^2r^2\right),\vspace{2.3mm}\\
		\dfrac{{\partial}^3 P_1(a,r)}{\partial a^3}=-12r^3(1+2ar)\leq 0\;\;\mbox{for all}\; r,\; a\in[0,1).
	\end{cases}
	\]
	We see that $\frac{{\partial}^2 P_1(a,r)}{\partial a^2}$ is a decreasing function of $a$, and hence we see that $ \frac{{\partial}^2 P_1(a,r)}{\partial a^2}\leq \frac{{\partial}^2 P_1(0,r)}{\partial a^2}= 0\;\;\mbox{for all}\; r\in [0,1)\;\; \mbox{and}\;\; a\in [0,1). $ Clearly, $ \frac{\partial P_1(a,r)}{\partial a}\geq \frac{\partial P_1(1,r)}{\partial a}=5r(1+r)\left({1}/{\sqrt{5}}+r\right)\left({1}/{\sqrt{5}}-r\right)\geq 0\;\;\mbox{for all}\;\; r\leq {1}/{\sqrt{5}} $ which shows that $P_1(a,r)$ is an increasing function $a$. Therefore, it is easy to see that $ 	P_1(a,r)\leq P_1(1,r)=-(1-r^2)^2\leq 0\;\;\mbox{for every }\;\;r\leq {1}/{\sqrt{5}}. $ This proves that $\mathcal{M}_{1, f}(r)\leq1$ for $r\leq 1/\sqrt{5}$ and $a\in[0,1).$ In order to show that the constant $r= 1/\sqrt{5}$ is sharp, we consider the function 
	\begin{align}\label{ee-3.12}
		\Phi_b(z)&:=\left(\dfrac{b-z}{1-bz}\right)I=A_0+\sum_{n=1}^{\infty}A_nz^n,\;\; z\in\mathbb{D},
	\end{align}
	where $A_0=bI,$ $A_n=-(1-b^2)b^{n-1}I$ and for some $b\in[0,1).$ With the help of the function $\Phi_b $, an easy computation shows that
	\begin{align*}
		\mathcal{M}_{1,\Phi_b}(r)&=||\Phi_b(-r)||+\sum_{n=1}^{\infty}||A_n||r^n+\left(\dfrac{1}{1+||A_0||}+\dfrac{r}{1-r}\right)\sum_{n=1}^{\infty}{||A_n||}^2r^{2n}\\&=\left(\dfrac{b+r}{1+br}\right)+\dfrac{(1-b^2)r}{1-br}+\dfrac{(1+br)(1-b^2)r^2}{(1+b)(1-r)(1-b^2r^2)}=1+\dfrac{(1-b)G(b,r)}{(1+br)(1-r)},
	\end{align*}
	where $G(b,r):=r(1+b)(1+br)-(1-r)^2.$ Clearly, $ \lim_{b\rightarrow 1^{-}} G(b,r)=-1+4r+r^2>0\;\;\mbox{for }\;\; r>1/\sqrt{5} $
	which turns out that $\mathcal{M}_{1,\Phi_b}(r)>1$ if $r>1/\sqrt{5}$ and as the limiting case $ b\rightarrow 1 $ suggests. This proves the sharpness of the constant $1/\sqrt{5}$  and the proof is complete. Again, we can easily get $\mathcal{M}_{2,f}\leq 1$ for $r\leq1/3$ and $1/3$ is the best possible,  since the proof of it follows on the similar as $\mathcal{M}_{1,f}$,  we omit the details and with this, proof of the lemma is completed.
\end{proof}
\begin{proof}[\bf Proof of lemma \ref{lem-4.3}.]
	For $f(z)=\sum_{n=0}^{\infty}A_nz^n\in\mathcal{B}(\mathbb{D},\mathcal{B}(\mathcal{H})), $ we obtain that 
	\begin{align}\label{ee-2.15}
		||f(z)-A_0||=||\sum_{n=0}^{\infty}A_nz^n-A_0||=||\sum_{n=01}^{\infty}A_nz^n||\leq\sum_{n=1}^{\infty}||A_n||r^n.
	\end{align} 
	In view of \eqref{ee-3.5} and \eqref{ee-2.15} with $a=||A_0||,$ we easily obtain
	\begin{align*}
		\mathcal{C}_{1, f}(r)&\leq a+2\dfrac{(1-a^2)r}{1-r}+\dfrac{(1+ar)(1-a^2)^2r^2}{(1+a)(1-r)(1-r^2)}=1+\dfrac{(1-a)P_4(a,r)}{(1-r)(1-r^2)},
	\end{align*}
	where $ P_4(a,r):=-(1-r)(1-r^2)+2r(1-r^2)(1+a)+r^2(1-a^2)(1+ra). $ By the similar arguments used in case of $ P_1(a, r) $, it can be easily shown that $P_4(a,r) $ is an increasing function of $a$ on the interval $[0,1)$, and hence $ P_4(a,r)\leq P_4(1,r)=(1-r^2)(5r-1) $
	which shows that $P_4(a,r)\leq 0$ for all $r\leq 1/5$ and $a\in [0,1)$. Consequently, the  desired inequality $ \mathcal{C}_{1, f}(r)\leq 1 $ holds.
	The sharpness of $ 1/5 $ can be easily shown by the similar arguments being used above for the function $\Phi_b$ given by \eqref{ee-3.12}, and hence we omit the details. Similarly, we can easily prove $ \mathcal{C}_{2, f}(r)\leq 1 $ holds for all $r\leq 1/3$ and $1/3$ is the best possible, since the proof of it follows on the similar as $\mathcal{C}_{1, f}(r)$,  we omit the details and with this, proof of the lemma is completed. 
\end{proof}

\begin{proof}[\bf Proof of Lemma \ref{lem-2.3}.]
	In view of \eqref{ee-3.5} with $a=||A_0||,$ we obtain that 
	\begin{align*}
		\mathcal{B}_p(r)&\leq a^p+\dfrac{(1-a^2)r}{1-r}+\left(\dfrac{1}{1+a}+\dfrac{r}{1-r}\right)\dfrac{(1-a^2)^2r^2}{(1-r^2)}:=Q(a),
	\end{align*}
	where $ Q(a)=a^p+k_1(1-a^2)+k_2(1-a)(1-a^2)+k_3(1-a^2)^2 $, $ 	k_1={r}/{(1-r)},\;\; k_2={r^2}/{(1-r^2)}\; \mbox{and}\; k_3={r^3}/((1-r)(1-r^2)) $. We claim that $Q(a)\leq 1$ for $r\leq 1/3$ and $a\in[0,1).$ A simple computation shows that
	\[
	\begin{cases}
		Q^{\prime}(a)=pa^{p-1}-2ak_1+k_2(3a^2-2a-1)+4k_3(a^3-a),\\
		Q^{\prime\prime}(a)=p(p-1)a^{p-2}-2k_1+k_2(6a-2)+4k_3(3a^2-1),\\ Q^{\prime\prime\prime}(a)=p(p-1)(p-2)a^{p-3}+6k_3+24k_3a.
	\end{cases} 
	\]
	We see that $Q^{\prime\prime\prime}(a)\geq 0$ for $a\in[0,1]$ and $p\in(0,1].$ This shows that $Q^{\prime\prime}(a)$ is an increasing function of $a,$ and hence $ Q^{\prime\prime}(a)\leq Q^{\prime\prime}(1)=-2k_1+4k_2+8k_3={2r(1+r)(3r-1)}/{((1-r)(1-r^2))}\leq 0,\;\;\mbox{for}\;\; r\leq{1/3}. $ This implies $Q^{\prime}(a)$ is a decreasing function of $a$ on $[0,1).$ A simple computation shows that $ Q^{\prime}(a)\geq Q^{\prime}(1)=1-2k_1={(1-3r)}/{(1-r)}\geq 0\;\;\mbox{for}\;\; r\leq {1}/{3}. $ This claim that $Q(a)$ is an increasing function of $a$ on $[0,1).$ Thus, $ Q(a)\leq Q(1)=1\;\;\mbox{for}\;\; r\leq {1}/{3}. $
	Therefore, the desired inequality \eqref{ee-4.4} established. 
	The sharpness of the constant $ 1/3 $ can be easily shown by the similar arguments being used above for the function $\Phi_b$ given by \eqref{ee-3.12}. Hence we omit the details and with this, proof of the lemma is completed.
\end{proof}

\begin{proof}[\bf Proof of Lemma \ref{lem-3.4}.]
	In view of \eqref{ee-3.5} and \eqref{e-2.8}  with $||A_0||=a,$ we obtain 
	\begin{align}\label{e-4.3}
		\dfrac{S_r}{\pi}\leq(1-a^2)^2\sum_{n=1}^{\infty}nr^{2n}=\dfrac{(1-a^2)^2r^2}{(1-r^2)^2}.
	\end{align}
	With the help of \eqref{ee-3.5}, \eqref{ee-3.6} and  \eqref{e-4.3}, an easy computation shows that
	\begin{align*}
		\mathcal{D}_{1, f}(r)\leq \left(\dfrac{a+r}{1+ra}\right)+\dfrac{(1-a^2)r}{1-r}+\sum_{m=1}^{\infty}c_m(1-a^2)^2\left(\dfrac{r}{1-r^2}\right)^{2m}:=\mathcal{M}^*(r).
	\end{align*}
	Clearly, $\mathcal{M}^*(r)$ is an increasing function of $r$, and hence we have $\mathcal{M}^*(r)\leq \mathcal{M}^*(\sqrt{5}-2)$ for $r\leq \sqrt{5}-2.$ A simple computation shows that 
	\begin{align*}
		&\mathcal{M}^*(\sqrt{5}-2)\\&\leq1-(1-a)\left(\dfrac{\left((3-\sqrt{5})^2-(\sqrt{5}-2)(1+a)(1+(\sqrt{5}-2)a)\right)}{(3-\sqrt{5})(1+(\sqrt{5}-2)a)}-\sum_{m=1}^{l}c_m2^{(1-4m)}\right)\\&\leq 1,
	\end{align*}
	if, and only, if
	\begin{align*}
		\sum_{m=1}^{l}c_m2^{(1-4m)} \leq\dfrac{\left((3-\sqrt{5})^2-(\sqrt{5}-2)(1+a)(1+(\sqrt{5}-2)a)\right)}{(3-\sqrt{5})(1+(\sqrt{5}-2)a)}:=\mathcal{V}(a).
	\end{align*}
	Since $\mathcal{V}(a)$ is a decreasing function of $a\in (0,1]$ and the maximal value of $\mathcal{V}(a)$ attain at $a=0$. Therefore, the maximal value is $(13-\sqrt{5})/4.$ Thus, the desired inequality \eqref{eee-3.8} is holds for $r\leq \sqrt{5}-2$ if  
	\begin{align*}
		\sum_{m=1}^{l}c_m2^{(1-4m)}\leq\dfrac{13-5\sqrt{5}}{4}.
	\end{align*}\vspace{1.2mm} 
	
	In order to show the sharpness of the constant $r_0=\sqrt{5}-2,$ we consider the function $\Phi_b$ given by \eqref{ee-3.12}. A simple computation now shows that $\mathcal{D}_{1, \Phi_b}(-r)=1-(1-b)\mathcal{U}(b,r),$ 
	where  
	\begin{align*}
		\mathcal{U}(b,r)=\dfrac{(1-r)}{1+br}-\dfrac{(1+b)r}{1-br}-\sum_{m=1}^{l}(1-b^2)^{2m-1}(1+b).
	\end{align*}
	A simple computation shows that $\lim_{b\rightarrow 1^-}\mathcal{U}(b,\sqrt{5}-2)=0.$ Therefore $\mathcal{U}(b,r)<0,$ when $r>\sqrt{5}-2$ and $b\rightarrow1^-.$ Hence $\mathcal{D}_{1, \Phi_b}(r)=1-(1-b)\mathcal{U}(b,r)>1$ for $r>\sqrt{5}-2,$ which shows that $\sqrt{5}-2$ is the best possible and with this, proof of the lemma is completed. 
\end{proof}
\begin{proof}[\bf Proof of Lemma \ref{BS-lem-3.5}.]
	Suppose that $f\in\mathcal{B}(\mathbb{D},\mathcal{B}(\mathcal{H}))$ with the expansion $f(z)=\sum_{n=0}^{\infty}A_nz^n$ in $\mathbb{D}$ such that $A_n\in\mathcal{B}(\mathcal{H})$ for all  $ n\in\mathbb{N}_0 $ and $f(0)=A_0I,$ $A_0\in\mathbb{D}.$ Since $||f(z)||\leq1$ in $\mathbb{D},$
	\begin{align*}
		f(z)\prec\varphi(z)=\left(\frac{A_0-z}{1-A_0z}\right)I=A_0I-\sum_{k=1}^{\infty}(1-A_0^2)A_0^{k-1}Iz^k.
	\end{align*}
	By the properties of subordination, we obtain
	\begin{align}\label{BS-eq-4.4}
		\dfrac{S_r}{\pi}=\sum_{k=1}^{\infty} k||A_k||^2r^{2k}\leq \sum_{k=1}^{\infty} k\left((1-||A_0||^2)||A_0||^{k-1}\right)^2r^{2k}=\frac{r^2(1-||A_0||^2)^2}{(1-||A_0||^2r^2)^2} 
	\end{align}
for $ r\leq {1}/{\sqrt{2}} $. Thus, the lemma's proof is now concluded.
\end{proof}
\begin{proof}[\bf Proof of Lemma\ref{BS-lem-3.6}.]
	Suppose that $ ||A_0||=a $. In view of \eqref{ee-3.5}, \eqref{BS-eq-4.4} and through a basic computation, it can be shown that
	\begin{align*}
		\mathcal{N}_f^1(r)&=\sum_{n=0}^{\infty}||A_n||r^n+\left(\dfrac{1}{1+||A_0||}+\dfrac{r}{1-r}\right)\sum_{n=1}^{\infty}||A_n||^2r^{2n}+\dfrac{8}{9}\left(\dfrac{S_r}{\pi}\right)\\&\leq a+\dfrac{(1-a^2)r}{1-r}+\dfrac{1+ar}{(1+a)(1-r)}\dfrac{(1-a^2)^2r^2}{(1-r^2)}+\frac{8}{9}\dfrac{(1-a^2)^2r^2}{(1-a^2r^2)^2}\\&=1-(1-a)H_1(r),
	\end{align*}
	where 
	\begin{align*}
		H_1(r):=1-\frac{1+a}{\frac{1}{r}-1}-\frac{(\frac{1}{r}+a)(1-a^2)}{(\frac{1}{r}-1)(\frac{1}{r^2}-1)}-\frac{8}{9}\frac{\frac{1}{r^2}(1-a^2)(1+a)}{(\frac{1}{r^2}-a^2)^2}.
	\end{align*}
	By an elementary computation, it can be shown that the function $ H_1 $ is an increasing on $ (0, 1) $, hence we have $ H_1(r)\leq H_1(1/3) $. Certainly, to show that $ \mathcal{A}_f(r)\leq 1 $, it is enough to show that $ H_1(1/3)\leq 1 $. In fact, we see that 
	\begin{align*}
		H_1\left(\frac{1}{3}\right)=1-\frac{(1-a)^2}{16(9-a^2)^2}(277-857a+281a^2+371a^3-49a^4-27a^5+3a^6+a^7)\leq 1
	\end{align*}
	if, and only if, $ 0\leq a=||A_0||\leq 0.402964$. Consequently, we see that $ H_1(r)\leq 1 $ for $ r\leq 1/3 $.\vspace{1.2mm}
	It remains to show that $ 1/3 $ is best possible, we consider the function $ \Phi_b $ given by \eqref{ee-3.12}. For this $ \Phi_b, $ we see that
	\begin{align*}
		\mathcal{N}^1_{\Phi_b}(r)=&b+\frac{(1-b^2)r}{1-rb}+\frac{1+br}{(1+b)(1-r)}\frac{(1-b^2)^2r^2}{1-b^2r^2}+\lambda\frac{(1-b^2)^2r^2}{(1-b^2r^2)^2}.
	\end{align*}
	For $ r=1/3 $, $ \mathcal{N}^1_{{\Phi_b}}(1/3) $ becomes
	\begin{align*}
		\mathcal{N}^1_{\Phi_b}\left(\frac{1}{3}\right)=1+\frac{(1-b)^2}{2(9-b^2)^2}\Psi_{\lambda}(b),
	\end{align*}
	where 
	\begin{align*}
		\Psi_{\lambda}(b):=8(9\lambda-8)-8(9\lambda+4)(1-b)+6(3\lambda+2)(1-b)^2+4(1-b)^3-(1-b)^4.
	\end{align*}
	It is easy to see that $ \mathcal{N}^1_{\Phi_b}\left({1}/{3}\right)>1 $ in case $ \lambda>8/9 $ and $ b\rightarrow 1^{-} $. This proves that $1/3$ and $8/9$ are best possible.\vspace{1.2mm} 
	
	Moreover, considering \eqref{ee-3.5} and \eqref{BS-eq-4.4}, and employing a straightforward calculation, it becomes evident that
	\begin{align*}
		\mathcal{N}^2_f(r)&\leq a^2+(1-a^2)\dfrac{r}{1-r}+\dfrac{(1+ar)}{(1+a)(1-r)}\dfrac{(1-a^2)^2r^2}{1-r^2}+\dfrac{9}{8}\dfrac{r^2(1-a^2)^2}{(1-a^2r^2)^2}\\&=1-(1-a^2)\left(1-\dfrac{r}{1-r}-\dfrac{(1+ar)(1-a)r^2}{(1-r)(1-r^2)}-\dfrac{9r^2(1-a^2)}{8(1-a^2r^2)^2}\right).
	\end{align*}
	By the similar arguments, to show that $ \mathcal{N}^2_f(r)\leq 1 $, it is enough to show that $ \mathcal{N}^2_f(1/(3-a))\leq 1 $. In fact,
	an easy computation shows that
	\begin{align*}
		&\mathcal{N}^2_f\left(\frac{1}{3-a}\right)\\&\leq 1-\frac{(1-a^2)(1-a)\left(216-780a+876a^2-419a^3+95a^4-13a^5+a^6\right)}{8((3-2a)^2)(2-a)(a^2-6a+8)}\leq 1
	\end{align*}
	if, and only if, $ 0\leq a<0.489758. $\vspace{1.2mm}
	
	It remains to show that $ 1/(3-a) $ and $9/8$ are best possible, we consider the function $ \Phi_a $ with some $a\in (0,1)$ given by \eqref{ee-3.12}. For this $ \Phi_a, $ we see that
	\begin{align*}
		\mathcal{N}^2_{\phi_a}\left(\frac{1}{3-a}\right)=1+\frac{(1-a)^2(1+a)}{9(3-2a)(2-a)}\mathcal{N}^*_{\phi_a}(a,\lambda),
	\end{align*}
	where 
	\begin{align*}
		\mathcal{N}^*_{\phi_a}(a,\lambda):&=(8\lambda-9)+12(\lambda-3)(1-a)+2(\lambda-18)(1-a)^2-3\lambda(1-a)^3\\&\quad-\lambda(1-a)^4.
	\end{align*}
	It is easy to see that 
	\begin{align*}
		\mathcal{N}^2_{\phi_a}\left(\frac{1}{3-a}\right)>1\; \mbox{for}\; \lambda>9/8\; \mbox{and}\; a\rightarrow 1^{-}.
	\end{align*}
	With this, the proof is completed. 
\end{proof}
\begin{proof}[\bf Proof of Theorem \ref{thm-3.1}.]
	To obtain the inequality $ \mathcal{I}_{1, f}(r)\leq 1 $, we convert the multidimensional power series \eqref{e-2.1} into the power series of one complex variable, and then we will use the Lemma \ref{lem-33.11}.  Let $ 	L=\{z=(z_1,z_2,\dots,z_n): z_j=b_jh, j=1,2,\dots,n;h\in\mathbb{C}\} $	be a complex line. Then, in each section of the domain $Q$ by the line $L,$ the series \eqref{e-2.1} covert into the following power series of complex variable $h:$
	\begin{align*}
		f(bh)=\sum_{k=0}^{\infty}P_k(b)h^k=f(0)+\sum_{k=1}^{\infty}P_k(b)h^k.
	\end{align*}
	Since $||f(bh)||<1$ for $h\in\mathbb{D}$ and $f(0)=b_0I,$ $b_0\in\mathbb{C},$ by used of the Lemma \ref{lem-33.11}, we obtain 
	\begin{align}\label{ee-3.13}
		||f(bh)||+\sum_{k=1}^{\infty}{||P_k(b)h^k||}+\left(\dfrac{1}{1+||f(0)||}+\dfrac{r}{1-r}\right)\sum_{k=1}^{\infty}{||P_k(b)h^k||}^2\leq1 
	\end{align} 
	for $z$ in the section $L\cap(1/\sqrt{5})Q.$ Since $L$ is an arbitrary complex line passing through the origin, the inequality \eqref{ee-3.13} is just same as $ \mathcal{I}_{1, f}(r)\leq 1 $.\vspace{1.2mm} 
	
	For the sharpness of the constant  $1/\sqrt{5},$ let the domain $Q$ be convex. Then $Q$ is an intersection of half spaces $ Q=\cap_{b\in J}\{z=(z_1,z_2,\dots,z_n):Re(b_1z_1+\dots+b_nz_n)<1\} $ for some $J.$ Since $Q$ is a circular domain, we obtain $ 	Q=\cap_{b\in J}\{z=(z_1,z_2,\dots,z_n):|b_1z_1+\dots+b_nz_n|<1\}. $
	To show that the constant $1/\sqrt{5}$ is the best possible, it is suffices to show that $1/\sqrt{5}$ cannot be improved for each domain $ 	R_b:=\{z=(z_1,z_2,\dots,z_n):|b_1z_1+\dots+b_nz_n|<1\} $ In view of the Lemma \ref{thm-3.1}, for some $b\in[0,1),$  there exist a function $\Phi_b:\mathbb{D}\rightarrow\mathcal{B}(\mathcal{H})$ defined by \eqref{ee-3.12} with $||\Phi_b(z)||<1$ for $z\in\mathbb{D},$  but for every $|z|=r>1/\sqrt{5},$ it is easy to see that $\mathcal{M}_{1, \Phi_b}(r)\leq1$ fails to hold in the disk $\mathbb{D}_r=\{z:|z|<r\}.$ On the other hand, consider the function $\omega:R_b \rightarrow \mathbb{D}$  defined by $\omega(z):=b_1z_1+\dots+b_nz_n.$ Thus, the function $f(z)=(\Phi_b\;\circ\; \omega)(z)$ gives the sharpness of the constant $1/\sqrt{5}$ for each domain  $R_b.$ \vspace{1.2mm}
	
	Now, by the similar argument, with the help of the Lemma \ref{lem-33.11} and the analogues proof of the inequality $ \mathcal{I}_{1, f}(r)\leq 1 $ and in view of $ \mathcal{M}_{2, f}(r)\leq 1 $, we can obtain the inequality $ \mathcal{I}_{2, f}(r)\leq 1 $ in the homothetic domain $(1/3)Q.$ To prove the constant $1/3$ is best possible when $Q$ is convex, in view of the analogues proof of the inequality $ \mathcal{I}_{1, f}(r)\leq 1 $, it is enough to show that the constant $1/3$ cannot be improved for each domain $ R_b $. Therefore, the function $f(z)=(\Phi_b\;\circ\; \omega)(z)$ gives the sharpness of the constant $1/3$  in each domain $R_b$. With this the proof is completed. 
\end{proof}

\begin{proof}[\bf Proof of Theorem \ref{thm-3.3}]
	In view of Lemma \ref{lem-4.3} and the analogues proof of Theorem \ref{thm-3.1}, we can obtain the inequality $ \mathcal{J}_{1, f}(r)\leq 1 $ in the homothetic domain $(1/5)Q,$ and inequality $ \mathcal{J}_{2, f}(r)\leq 1 $ for $z$ in the homothetic domain $(1/3)Q.$ To prove the constants $1/5$ and $1/3$ are the best possible when $Q$ is convex, in view of the analogues proof of Theorem \ref{thm-3.1}, it is enough to show that $1/5$ and $1/3$ cannot be improved for each domain  $R_b$ Thus, the function  $f(z)=(\Phi_b\;\circ\; \omega)(z)$ gives the sharpness of the constants $1/5$ and $1/3$ in each domain $R_b.$ This completes proof of the theorem.   
\end{proof}

\begin{proof}[\bf Proof of Theorem \ref{thm-3.2}]
	By means of \eqref{ee-4.4} of Lemma \ref{lem-2.3}, using the analogous proof of Theorem \ref{thm-3.1}, we can easily obtain the inequality \eqref{ee-33.5} in the homothetic domain $(1/3)Q.$ To prove the constant $1/3$ is best possible when $Q$ is convex,  in view of the analogous proof of the Theorem \ref{thm-3.1}, it is sufficient to show that the number $1/3$ cannot be improved for each of the domain $R_b$. Therefor, the function $f(z)=(\Phi_b\circ \omega)(z)$ gives the sharpness of the constant $1/3$ in each domain $R_b$. The proof is completed.  
\end{proof}	

\begin{proof}[\bf Proof of Theorem \ref{thh-3.4}]
	Using Lemma \ref{lem-3.4} and the analogues proof of the Theorem \ref{thm-3.1}, we can obtain the inequality  $ \mathcal{L}_{ f}(r)\leq 1 $  in the homothetic domain $(\sqrt{5}-2)Q ,$  where the coefficients of $G$ satisfy \eqref{eee-3.8}. Now we prove that if $Q$ is convex, then the constants $\sqrt{5}-2$ cannot be improved. The reasoning as in the proof of Theorem \ref{thm-3.1}  concludes the proof of the inequality $ \mathcal{L}_{ f}(r)\leq 1 $ and with this the proof of the theorem completes.    
\end{proof}	
\begin{proof}[\bf Proof of Theorem \ref{BS-thm-3.5}]
	By means of Lemma \ref{BS-lem-3.6} and utilizing the analogous proof presented in Theorem \ref{thm-3.1}, we can verify that $\mathcal{R}_f^1(r)\leq 1$ holds for $z$ in the homothetic domain $(1/3)Q$, as well as $\mathcal{R}_f^2(r)\leq 1$ for $z$ in the homothetic domain $(1/(3-||A_0||))Q$. Our next task is to prove that, considering the convexity of $Q$, the constants $1/3$, $8/9$, $9/8$, and $1/(3-||A_0||)$ cannot be improved. The same reasoning as in the proof of Theorem \ref{thm-3.1} concludes the proof of the inequalities $\mathcal{R}^1_{f}(r)\leq 1$ and $\mathcal{R}^2_{f}(r)\leq 1$. With this, the theorem's proof is complete.   
\end{proof}	

\end{document}